\magnification=1200
\overfullrule=0pt
\centerline {\bf Uniqueness properties of functionals with Lipschitzian derivative}\par
\bigskip
\bigskip
\centerline {BIAGIO RICCERI}\par
\bigskip
\bigskip
\bigskip
\bigskip
\noindent
{\bf 1. Introduction}\par
\bigskip

Let $X$ be a real Hilbert space and $J$ a $C^1$ functional on $X$. For
$x_0\in X$, $r>0$, set $S(x_0,r)=\{x\in X : \|x-x_0\|=r\}$.\par
\smallskip
Also on the basis of the beautiful theory developed and applied by Schechter and
Tintarev in [2], [3], [4] and [5], it is of particular interest to know when
the restriction of $J$ to $S(0,r)$ has a unique maximum.\par
\smallskip
The aim of the present paper is to offer a contribution along this direction.\par
\smallskip
We show that such a uniqueness property holds (for suitable $r$) provided
that $J'$ is Lipschitzian and $J'(0)\neq 0$. At the same time, we also
show that (for suitable $s$) the set $J^{-1}(s)$ has a unique element of 
minimal norm.\par
\smallskip
After proving the general result (Theorem 1), we present an application to a semilinear
Dirichlet problem involving a Lipschitzian nonlinearity (Theorem 2).\par
\bigskip
\noindent
{\bf 2. The main result}\par
\bigskip
With the usual convention $\inf \emptyset=+\infty$, our main result reads as follows:
\par
\medskip
THEOREM 1. - {\it Let $X$ be a real Hilbert space and let $J:X\to {\bf R}$
be a sequentially weakly upper semicontinuous $C^1$ functional, with Lipschitzian
derivative. Let $L$ be the Lipschitz constant of $J'$.\par
Then, for each $x_0\in X$ with $J'(x_0)\neq 0$, if we set
$$\alpha_0=\inf_{x\in M_{1\over L}}J(x)$$ and
$$\beta_0=\hbox {\rm dist}(x_0,M_{1\over L})\ ,$$
where $M_{1\over L}$ is the set of all global minima of the functional
$x\to {{1}\over {2}}\|x-x_0\|^2-{{1}\over {L}}J(x)$, we have
$\alpha_0>J(x_0)$,
$\beta_0>0$,
and the following properties hold:\par
\noindent
$(i)$\hskip 5pt 
 for every $r\in ]J(x_0),\alpha_0[$ there exists a unique
$y_r\in J^{-1}(r)$ such that
$$\|x_0-y_r\|=\hbox {\rm dist}(x_0,J^{-1}(r))\ ;$$ 
$(ii)$\hskip 5pt for every $r\in ]0,\beta_0[$ the restriction
of the functional $J$ to the set $S(x_0,r)$ has a unique
global maximum.}
\medskip
 The main tool used to get Theorem 1 is the following 
particular case of Theorem 3 of [1].\par
\medskip
THEOREM A. - 
{\it Let $X$ be a reflexive real Banach space,
$I\subseteq {\bf R}$ an interval 
and $\Psi : X\times I\to {\bf R}$ a function such
that $\Psi(x,\cdot)$ is concave and continuous for all $x\in X$, while
$\Psi(\cdot,\lambda)$ is sequentially weakly lower semicontinuous 
and coercive, with a unique local minimum for all $\lambda\in
\hbox {\rm int}(I)$.\par
Then, one has
$$\sup_{\lambda\in I}\inf_{x\in X}\Psi(x,\lambda)=
\inf_{x\in X}\sup_{\lambda\in I}\Psi(x,\lambda)\ .$$}\par
\medskip
We will also use the two propositions below.\par
\medskip
PROPOSITION 1. - {\it Let $Y$ be a nonempty set, 
$f, g:Y\to {\bf R}$ two functions, and $a, b$ two
real numbers, with $a<b$. Let $y_a$ 
 be a global minimum of the function $f-ag$ and
$y_b$ a global minimum of the function
$f-bg$.\par
Then, one has $g(y_a)\leq g(y_b)$. If either $y_a$
or $y_b$ is strict and $y_a\neq y_b$, then
$g(y_a)<g(y_b)$.}\par
\smallskip
PROOF. We have
$$f(y_a)-ag(y_a)\leq f(y_b)-ag(y_b)$$
as well as
$$f(y_b)-bg(y_b)\leq f(y_a)-bg(y_a)\ .$$
Summing, we get
$$-ag(y_a)-bg(y_b)\leq -ag(y_b)-bg(y_a)$$
and so
$$(b-a)g(y_a)\leq (b-a)g(y_b)$$
from which the first conclusion follows.
If either $y_a$ or $y_b$ is strict and $y_a\neq y_b$,
then one of the first two inequalities is strict and hence
so is the third one.
\hfill $\bigtriangleup$\par
\medskip
PROPOSITION 2. - {\it Let $Y$ be a real Hibert space
and let $\varphi:Y\to {\bf R}$ be a sequentially
weakly upper semicontinuous $C^1$ functional whose
derivative is a contraction.\par
Then, for every $y_0\in Y$, the functional
$y\to {{1}\over {2}}\|y-y_0\|^2-\varphi(y)$ is
coercive and has a unique local minimum.} \par
\smallskip
PROOF. Let $\nu$ be the Lipschitz constant of $\varphi'$.
So, $\nu<1$, by assumption. For each $y\in Y$, we have
$$\varphi(y)=\varphi(0)+\int_{0}^{1}\langle \varphi'(ty),y\rangle dt$$
and so
$$|\varphi(y)|\leq |\varphi(0)|+\int_{0}^{1}|\langle \varphi'(ty),y
\rangle |dt\leq |\varphi(0)|+\|y\|\int_{0}^{1}\|\varphi'(ty)\|dt\leq$$
$$\leq |\varphi(0)|+\|y\|\left ( \int_{0}^{1}\|\varphi'(ty)-\varphi'(0)\|dt
+\|\varphi'(0)\|\right )\leq |\varphi(0)|+{{\nu}\over {2}}\|y\|^2+
\|\varphi'(0)\|\|y\|\ .$$
 From this, we then get
$${{1}\over {2}}\|y-y_0\|^2-\varphi(y)\geq {{1-\nu}\over {2}}\|y\|^2-
(\|\varphi'(0)\|+\|y_0\|)\|y\|+{{1}\over {2}}\|y_0\|^2-|\varphi(0)|$$
and hence
$$\lim_{\|y\|\to +\infty}{{1}\over {2}}\|y-y_0\|^2-\varphi(y)=+\infty$$
which yields our first claim. Then, the functional
$y\to {{1}\over {2}}\|y-y_0\|^2-\varphi(y)$ has a global minimum, since
it is sequentially weakly lower semicontinuous. But the critical points of
this functional are exactly the fixed points of the operator $\varphi'+y_0$
that is a contraction. So, the functional has a unique local minimum (that is its
global minimum).\hfill $\bigtriangleup$
\medskip
{\it Proof of Theorem 1}. First, note that, for each $\gamma>L$, the operator
${{1}\over {\gamma}}J'$ is a contraction, and so,
by Proposition 2, the functional
$x\to {{1}\over {2}}\|x-x_0\|^2-{{1}\over {\gamma}}J'(x)$ has a unique global
minimum, say $x_{1\over \gamma}$. Fix $\gamma>L$. By Proposition 1, we have
 $J(x_0)\leq J(x_{1\over \gamma})$. We claim that $J(x_0)<J(x_{1\over \gamma})$.
Arguing by contradiction, assume that $J(x_0)=J(x_{1\over \gamma})$. Then,
by Proposition 1 again, we would have $x_{1\over \gamma}=x_0$. 
 Consequently, the derivative of the functional
$x\to {{1}\over {2}}\|x-x_0\|^2-{{1}\over {\gamma}}J'(x)$ would vanish at $x_0$,
that is $-{{1}\over {\gamma}}J'(x_0)=0$, against one of the hypotheses.
Then,
we have
$$J(x_0)<J(x_{1\over \gamma})\leq J(x)$$
for all $x\in M_{1\over L}$, and so $J(x_0)<\alpha_0$. Clearly,
$x_{1\over \gamma}$ is the global minimum
of the functional $x\to {{\gamma}\over {2}}\|x-x_0\|^2-J(x)$,
while any $z\in M_{1\over L}$ is a global minimum of the functional
$x\to {{L}\over {2}}\|x-x_0\|^2-J(x)$. Consequently, if we apply
Proposition 1 again
(with $f(x)=-J(x)$, $g(x)=-\|x-x_0\|^2$, $a={{L}\over {2}}$,
$b={{\gamma}\over {2}}$), for any $z\in M_{1\over L}$, we get
$$-\|z-x_0\|^2\leq -\|x_{1\over \gamma}-x_0\|^2\ ,$$
and so
$$\beta_{0}\geq \|x_{1\over \gamma}-x_0\|>0\ .$$
Now, to prove $(i)$, fix $r\in ]J(x_0),\alpha_0[$ and
consider the function $\Psi:X\times [0,{{1}\over {L}}]\to
{\bf R}$ defined by
$$\Psi(x,\lambda)={{1}\over {2}}\|x-x_0\|^2+\lambda (r-J(x))$$
for all $(x,\lambda)\in X\times [0,{{1}\over {L}}]$.
Taken Proposition 2 into account, it is clear that the function $\Psi$
satisfies all the assumptions of Theorem A. Consequently, we have
$$\sup_{\lambda\in [0,{{1}\over {L}}]}\inf_{x\in X}
\Psi(x,\lambda)=
\inf_{x\in X}\sup_{\lambda\in [0,{{1}\over {L}}]}
\Psi(x,\lambda)\ .$$
The functional $\sup_{\lambda\in [0,{{1}\over {L}}]}\Psi(\cdot,\lambda)$
is sequentially weakly lower semicontinuous and coercive, and so
there exists $x^*\in X$ such that
$$\sup_{\lambda\in [0,{{1}\over {L}}]}\Psi(x^*,\lambda)=
\inf_{x\in X}\sup_{\lambda\in [0,{{1}\over {L}}]}\Psi(x,\lambda)\ .$$
Also, the function $\inf_{x\in X}\Psi(x,\cdot)$ is upper semicontinuous,
and so there exists $\lambda^*\in [0,{{1}\over {L}}]$ such that
$$\inf_{x\in X}\Psi(x,\lambda^*)=\sup_{\lambda\in [0,{{1}\over {L}}]}
\inf_{x\in X}\Psi(x,\lambda)\ .$$
Hence, from this it follows that
$${{1}\over {2}}\|x^*-x_0\|^2+\lambda^* (r-J(x^*))=
\inf_{x\in X}{{1}\over {2}}\|x-x_0\|^2+\lambda^* (r-J(x))=
\sup_{\lambda\in [0,{{1}\over {L}}]}
{{1}\over {2}}\|x^*-x_0\|^2+\lambda (r-J(x^*))\ .$$
We claim that $J(x^*)=r$. Indeed, if it were $J(x^*)<r$, then we
would have $\lambda^*={{1}\over {L}}$, and so $x^*\in M_{1\over L}$,
against the fact that $r<\alpha_0$. If it were $J(x^*)>r$, then
we would have $\lambda^*=0$, and so $x^*=x_0$, against the fact that
$J(x_0)< r$. We then have
$${{1}\over {2}}\|x^*-x_0\|^2=
\inf_{x\in X}{{1}\over {2}}\|x-x_0\|^2+\lambda^*(r-J(x))\ .$$
This implies, on one hand, that $\lambda^*<{{1}\over {L}}$ (since
$r<\alpha_0$) and, on the other hand, that each global minimum
(and $x^*$ is so) of the restriction to $J^{-1}(r)$ of the functional
$x\to {{1}\over {2}}\|x-x_0\|^2$ 
is a global minimum
in $X$ of the functional
$x\to {{1}\over {2}}\|x-x_0\|^2-\lambda^*J(x)$.
 But this functional (just because
$\lambda^*<{{1}\over {L}}$) has a unique global minimum, and so
$(i)$ follows. Let us now prove $(ii)$. To this end, fix
$r\in ]0,\beta_{0}[$ and consider
the function $\Phi:X\times [L,+\infty[\to {\bf R}$ defined by
$$\Phi(x,\lambda)={{\lambda}\over {2}}(\|x-x_0\|^2-r^2)-J(x)$$
for all $(x,\lambda)\in X\times [L,+\infty[$. Applying Theorem A, we get
$$\sup_{\lambda\in [L,+\infty[}\inf_{x\in X}\Phi(x,\lambda)=
\inf_{x\in X}\sup_{\lambda\in [L,+\infty[}
\Phi(x,\lambda)\ .$$
Arguing as before (note, in particular,
 that $\lim_{\lambda\to +\infty}\inf_{x\in X}\Phi(x,\lambda)=-\infty$),
we get $\hat x\in X$ and $\hat \lambda\in [L,+\infty[$ such that
$$\sup_{\lambda\in [L,+\infty[}\Phi(\hat x,\lambda)=\inf_{x\in X}\sup_{\lambda\in [L,+\infty[}
\Phi(x,\lambda)$$
and
$$\inf_{x\in X}\Phi(x,\hat \lambda)=\sup_{\lambda\in [L,+\infty[}\inf_{x\in X}\Phi(x,\lambda)\ .$$
So that
$${{\hat \lambda}\over {2}}(\|\hat x-x_0\|^2-r^2)-J(\hat x)=
\inf_{x\in X}{{\hat \lambda}\over {2}}(\|x-x_0\|^2-r^2)-J(x)=
\sup_{\lambda\in [L,+\infty[}{{\lambda}\over {2}}(\|\hat x-x_0\|^2-r^2)-J(\hat x)\ .$$
 From this it follows at once that $\|\hat x-x_0\|^2\leq r^2$. But, if it were
$\|\hat x-x_0\|^2<r^2$ we would have $\hat \lambda=L$. This, in turn,
would imply that $\hat x\in M_{1\over L}$, against the fact that
$r<\beta_0$. Hence, we have $\|\hat x-x_0\|^2=r^2$. Consequently
$$-{{1}\over {\hat \lambda}}J(\hat x)=
\inf_{x\in X}{{1}\over {2}}(\|x-x_0\|^2-r^2)-{{1}\over {\hat \lambda}}J(x)\ .$$
This implies, on one hand, that $\hat \lambda>L$ (since $r<\beta_{0}$) and, on the
other hand, that each global maximum (and $\hat x$ is so)
of the restriction of the functional $J$
 to the set $S(x_0,r)$
 is a global mimimum in $X$ of the functional
$x\to {{1}\over {2}}\|x-x_0\|^2-{{1}\over {\hat \lambda}}J(x)$. Since $\hat \lambda>L$,
this functional has a unique global minimum, and so $(ii)$ follows.\hfill $\bigtriangleup$
\par
\medskip
REMARK 1. - It is clear from the proof that the assumption $J'(x_0)\neq
0$ has been used to prove $\alpha_0>J(x_0)$ and $\beta_0>0$, while
it has no role in showing $(i)$ and $(ii)$.
However, when $J'(x_0)=0$, it can happen that $\alpha_0=J(x_0)$,
$\beta_0=0$, with $(i)$  (resp. $(ii)$) holding for no $r>\alpha_0$
(resp. for no $r>0$). To see this, take, for instance, $X={\bf R}$,
$J(x)={{1}\over {2}}x^2$, $x_0=0$.  \par
\bigskip
\noindent
{\bf 3. An application}\par
\bigskip
  From now on, $\Omega$ is an open, bounded and connected subset of
${\bf R}^n$ with sufficiently smooth boundary, and
 $X$ denotes the space $W^{1,2}_{0}(\Omega)$, with the usual norm
$$\|u\|=\left ( \int_{\Omega}|\nabla u(x)|^{2}dx\right ) ^{1\over 2}\ .$$
Moreover, $f:{\bf R}\to {\bf R}$ is a Lipschitzian function,
with Lipschitz constant $\mu$.\par
\smallskip
Let $\lambda\in {\bf R}$. As usual, a classical solution of the
problem
$$\cases {-\Delta u=
\lambda f(u)
 & in
$\Omega$\cr & \cr u_{|\partial \Omega}=0 \cr} \eqno{(P_\lambda)}$$
is any $u\in C^2(\Omega)\cap C^0(\overline {\Omega})$, zero
on $\partial \Omega$, which satisfies the equation pointwise in $\Omega$.
\par
\smallskip
For each $u\in X$,  put
$$J(u)=
\int_{\Omega}\left ( \int_{0}^{u(x)}f(\xi)d\xi\right ) dx\ .$$
\smallskip
  By classical results, 
the functional $J$ is  continuously
G\^ateaux differentiable and sequentially weakly continuous in
$X$, and one has
$$J'(u)(v)=\int_{\Omega}f(u(x))v(x)dx$$
for all $u, v\in X$. 
Moreover, by a standard regularity result,
 the critical points in $X$ of the functional
$u\to {{1}\over {2}}\|u\|^2-\lambda J(u)$ are exactly
the classical solutions of problem $(P_\lambda)$.
\medskip
 Denote by $\lambda_{1}$ the first eigenvalue of the problem
$$\cases {-\Delta u=\lambda u & in
$\Omega$\cr & \cr u_{|\partial \Omega}=0 \ .\cr}$$
Recall that
$\|u\|_{L^{2}(\Omega)}\leq \lambda_{1}^{-{{1}\over {2}}}\|u\|$
for all $u\in X$. \par
\smallskip
We are now in a position to state the following\par
\medskip
THEOREM 2. - {\it Assume that $f(0)\neq 0$.
For each $r>0$, put
$$\gamma(r)=\sup_{\|u\|^2=r}J(u)\ .$$
Further, put
$$\delta_0=\inf_{u\in M}\|u\|^2$$
where $M$ is the set of all global minima in $X$ of
the functional $u\to {{1}\over {2}}\|u\|^2-{{\lambda_1}\over {\mu}}J(u)$.\par
Then, $\delta_0>0$, the function $\gamma$ is $C^1$ and $\gamma'$ is positive in $]0,\delta_0[$
and there exists a continuous function $\varphi:]0,\delta_0[\to X$ such that,
for each $r\in ]0,\delta_0[$, $\varphi(r)$ is a classical solution of the problem
$$\cases {-\Delta u=
{{1}\over {2\gamma'(r)}}f(u)
 & in
$\Omega$\cr & \cr u_{|\partial \Omega}=0\ \cr} $$
satisfying $\|\varphi(r)\|^2=r$ and $J(\varphi(r))=\gamma(r)$.}\par
\smallskip
PROOF.  
Fix $u, v, w\in X$, with $\|w\|=1$.
 We have
$$|J'(u)(w)-J'(v)(w)|\leq 
\int_{\Omega}|f(u(x))-f(v(x))||w(x)|dx\leq
\mu\|u-v\|_{L^{2}(\Omega)}\|w\|_{L^{2}(\Omega)}\leq
{{\mu}\over {\lambda_1}}\|u-v\|\ ,$$
and hence
$$\|J'(u)-J'(v)\|\leq {{\mu}\over {\lambda_1}}\|u-v\|\ .$$
That is, $J'$ is Lipschitzian in $X$, with Lipschitz constant
${{\mu}\over {\lambda_1}}$. Moreover, since $f(0)\neq 0$,
we have $J'(u)\neq 0$ for all $u\in X$.
Then, thanks to Theorem 1, for each
$r\in ]0,\delta_0[$, the restriction
of the functional $J$ to the sphere $S(0,\sqrt r)$ has a unique
maximum. At this point, taken into account that $\gamma(r)>0$ for all
$r>0$,
the conclusion follows directly from Lemma 2.1 and Corollary 2.13 of [2].
\hfill $\bigtriangleup$

\bigskip
\bigskip
\centerline {\bf References}\par
\bigskip
\bigskip
\noindent
[1]\hskip 5pt B. RICCERI, {\it Minimax theorems for limits of parametrized
functions having at most one local minimum lying in a certain set}, preprint.
\par
\smallskip
\noindent
[2]\hskip 5pt M. SCHECHTER and K. TINTAREV, {\it Spherical maxima in
Hilbert space and semilinear elliptic eigenvalue problems}, Differential
Integral Equations, {\bf 3} (1990), 889-899.\par
\smallskip
\noindent
[3]\hskip 5pt M. SCHECHTER and K. TINTAREV, {\it Points of shperical
maxima and solvability of semilinear elliptic equations}, Canad. J.
Math., {\bf 43} (1991), 825-831.\par
\smallskip
\noindent
[4]\hskip 5pt M. SCHECHTER and K. TINTAREV, {\it Eigenvalues for semilinear
boundary value problems}, Arch. Rational Mech. Anal., {\bf 113} (1991),
197-208.\par
\smallskip
\noindent
[5]\hskip 5pt M. SCHECHTER and K. TINTAREV, {\it Families of 'first
eigenfunctions' for semilinear elliptic eigenvalue problems}, Duke
Math. J., {\bf 62} (1991), 453-465.\par

\bigskip
\bigskip
\bigskip
\bigskip
Department of Mathematics\par
University of Catania\par
Viale A. Doria 6\par
95125 Catania\par
Italy\par
{\it e-mail address}: ricceri@dmi.unict.it

\bye